\documentclass[11pt,a4paper,twoside]{amsart}
\usepackage{amssymb,latexsym,amsmath,graphics}
\usepackage{graphicx}
\usepackage{caption}
\usepackage{subcaption}
\usepackage{enumerate}

\numberwithin{equation}{section}

\newcommand{\newtext}{}
 
\newcommand{\mntext}{}

\newtheorem {theorem}{Theorem}
\newtheorem {lemma}[theorem]{Lemma}
\newtheorem {proposition}[theorem]{Proposition}
\newtheorem {definition}[theorem]{Definition}

\theoremstyle{definition}

\newcommand{\argmin}{\mathop{\mathrm{arg\,min}}}
\newcommand{\Glim}{\Gamma\!\hbox{-}{\rm lim}}

\newcommand{\finu}{\mathbf{u}}
\newcommand{\finm}{\mathbf{m}}
\newcommand{\finA}{\mathbf{A}}
\newcommand{\finI}{\mathbf{I}}
\newcommand{\finL}{\mathbf{L}}
\newcommand{\finD}{\mathbf{D}}
\newcommand{\fineps}{\mathbf{e}}

\def \ba {\begin {eqnarray*} }
\def \ea {\end {eqnarray*} }
\def \beq {\begin {eqnarray}}
\def \eeq {\end {eqnarray}}

\newcommand{\T}{{\mathbb T}}
\newcommand{\R}{{\mathbb R}}
\newcommand{\E}{{\mathbb E}}
\newcommand{\N}{{\mathbb N}}
\newcommand{\Z}{{\mathbb Z}}

\newcommand{\ra}{\rightarrow}

\newcommand{\torus}{\T}

\newcommand{\beqno}{\begin{eqnarray*}}
\newcommand{\eeqno}{\end{eqnarray*}}

\newcommand{\ind}{\text{index} }
\newcommand{\Ker}{\text{Ker}\ }
\newcommand{\coker}{\text{Coker}\ }

\newcommand\opensubset{\mathrel{\ooalign{$\subset$\cr
\hidewidth\hbox{$\circ\mkern.5mu$}\cr}}}

\newcommand{\beqla}[1] {\begin {eqnarray}\label{#1}}

\def \eps {\varepsilon}

\def \beq {\begin {eqnarray}}
\def \eeq {\end {eqnarray}}
\def \ba {\begin {eqnarray*}}
\def \ea {\end  {eqnarray*}}

\def \b {\beta}
\def \a {\alpha}

\def \bra{\langle}
\def \cet{\rangle}

\def \noise {\varepsilon}

\def\tilde{\widetilde}

\title[Regularized inversion and  white Gaussian noise]{Analysis of regularized inversion of data \\ corrupted by white Gaussian noise}

\author{Hanne Kekkonen, Matti Lassas and Samuli Siltanen}
\address{Department of Mathematics and Statistics
P.O. Box 68,
00014 University of Helsinki
FINLAND}

\begin{document}

\maketitle

\begin{abstract}
{\newtext 
Tikhonov regularization is studied in the case of linear pseudodifferential operator as the forward map and additive white Gaussian noise as the measurement error. The measurement model for an unknown function $u(x)$ is 
\begin{eqnarray*}
m(x) = Au(x) + \delta\hspace{.2mm}\varepsilon(x),
\end{eqnarray*}
where $\delta>0$ is the noise magnitude. If $\varepsilon$ was an $L^2$-function, 
 Tikhonov regularization gives an estimate 
\begin{eqnarray*}
T_\alpha(m) = \argmin_{u\in H^r}\big\{\|A u-m\|_{L^2}^2+ \alpha\|u\|_{H^r}^2 \big\}
\end{eqnarray*}
for $u$ where $\alpha=\alpha(\delta)$ is the regularization parameter. Here penalization of the Sobolev norm $ \|u\|_{H^r}$ covers the cases of standard Tikhonov regularization ($r=0$) and first derivative penalty ($r=1$).

Realizations of white Gaussian noise are almost never in $L^2$, but do belong to $H^s$ with probability one if $s<0$ is small enough. A  modification of Tikhonov regularization theory is presented, covering the case of white Gaussian measurement noise. Furthermore, the convergence of regularized reconstructions to the correct solution as $\delta\rightarrow 0$ is proven in appropriate function spaces  using microlocal analysis. The convergence of
the related finite-dimensional problems to the infinite-dimensional problem is also analysed. }
%

\end{abstract}

 

\tableofcontents
  
\section{Introduction}\label{intro}  

\subsection{Discrete and continuous regularization}
\noindent 
Consider the following continuous model for indirect measurements:
\begin{equation}\label{contmodel1}
  m = Au + \mbox{noise},
\end{equation}
where the data $m$ and the quantity of interest $u$ are real-valued functions of $d$ real variables and $A$ is a bounded linear operator. A large class of practical measurements can be modelled by operators $A$ arising from partial differential equations of mathematical physics. We focus on ill-posed inverse problems where $A$ does not have a continuous inverse.

Physical measurement devices produce a discrete data vector $\finm\in\R^k$, which we model by adding a linear operator $P_k$ to (\ref{contmodel1}):
\begin{equation}\label{measprojectionPk}
  \finm := P_k(Au) + P_k(\mbox{noise}).
\end{equation}
Furthermore, practical solution of the inverse problem calls for a discrete representation of the unknown $u$. This can be done using some computationally feasible approximation of the form $\finu=T_n u\in\R^n$, for example Fourier series truncated to $n$ terms. 
The practical inverse problem is now 
\begin{equation}\label{discreteinverseproblem}
\mbox{\em given $\finm$, compute a noise-robust approximation to $\finu$}. 
\end{equation}

We study the most common computational appoach to (\ref{discreteinverseproblem}), namely classical Tikhonov regularization defined by
\begin{equation}\label{discrTikh1} 
  T_\alpha(\finm):=\argmin_{\finu\in\R^n} \left\{\| \finA\finu - \finm\|_2^2 + \alpha\|\finL\finu\|_2^2\right\}.
\end{equation}
Here $\finA=P_k A T_n$ is a $k\times n$ matrix approximation to the operator $A$, and $0<\alpha<\infty$ is the {\em regularization parameter}. The matrix $\finL$ is  used to introduce {\em a priori} information to the inversion. For example,
\begin{itemize}
\item[(a)] $\finL=\finI$, the identity matrix, models the {\em a priori} information that $u$ is not very large in norm.
\item[(b)] $\finL=\finI+\finD$, where $\finD$ is a finite-difference first-order derivative matrix,  models the {\em a priori} information that $u$ is continuously differentiable and $u$ or its derivative are not very large in square norm.
\end{itemize}
Our aim is to provide new analytic insight to the relationship between the continuous model  (\ref{contmodel1}) and practical inversion based on (\ref{discrTikh1}) in the case of {\newtext Gaussian noise.}
 
Note that the reconstruction $T_\alpha(\finm)$ given by (\ref{discrTikh1}) depends on both $k$ and $n$. Practical computational inversion may involve modifying both of them: updating the measurement device changes the number $k$ of data points, and refining the computational grid in the hope of extra accuracy will increase $n$. Furthermore, sometimes the most efficient numerical algorithm is based on a multigrid strategy, involving the computation of $T_\alpha(\finm)$ with several different $n$. 

Since there is a common continuous inverse problem behind the discrete model, it is desirable that the reconstruction $T_\alpha(\finm)$ converges to a meaningful limit as $k,n\rightarrow\infty$. Such convergence would also ensure that the dependency of $T_\alpha(\finm)$ on $k$ and $n$ is stable, at least for large enough values. Therefore, we discuss a continuous version of  (\ref{discrTikh1}) based directly on the ideal model (\ref{contmodel1}). 

{\newtext Under certain assumptions (including that $m$ should be an $L^2$-function) the finite-dimensional  problem (\ref{discrTikh1}) $\Gamma$-converges as $n,k\to \infty$ to the following infinite-dimensional minimization problem in a Sobolev space $H^r$:
\begin{equation}\label{contTikh00}
 \argmin_{u\in H^r}\big\{\|m-A u\|_{L^2}^2 + \alpha\|u\|_{H^r}^2 \big\}.
\end{equation}
See  Section \ref{sec:Gamma convergence}  below for a proof.
In (\ref{contTikh00}) the case}
$r=0$ corresponds to (a) and $r=1$ corresponds, roughly, to (b) above. 
However, formula (\ref{contTikh00}) only makes sense if the noise in (\ref{contmodel1}) is square integrable. This brings us to the main topic of the paper: noise modeling.
 
\subsection{Properties of white noise}
Next we will give the definitions for the discrete and continuous white noise and describe the 'white noise paradox' arising from the infinite $L^2$-norm of the natural limit of white Gaussian noise in $\R^k$ when $k\to\infty$.

We model the $k$-dimensional noise in (\ref{measprojectionPk}) as $P_k(\mbox{noise})=\delta\hspace{.2mm}\fineps$, where $\delta>0$ plays the role of noise amplitude. The vector $\fineps\in\R^k$ is a realization of a $\R^k$-valued Gaussian random variable ${\bf E}={\bf E}^{(k)}$ having  mean zero and unit variance: ${\bf E}^{(k)}\sim N(0,I)$. In terms of a probability density function we have
\begin{equation}\label{gauss}
\pi_{{\bf E}^{(k)}}(E)=c \, \hbox{exp}\bigg(-\frac {1}{2} \|E\|_2^2\bigg),\quad E\in \R^k,\quad \|E\|_2=(\sum_{j=1}^k E_j^2)^{1/2}.
\end{equation}
The appearance of $\|\,\cdotp\|_2$ in (\ref{gauss}) is the reason why  square norm is used in the data fidelity term $\| \finA\finu - \finm\|_2^2$ of (\ref{discrTikh1}).
The above noise model is appropriate for example for photon counting under high radiation intensity,
see e.g. \cite{dentalpaper2,dentalpaper1}.  

 {\newtext 
Let us relate the above to the continuous model (\ref{contmodel1}). We take $u(x)$ and $m(x)$ to be functions defined on a closed, compact $d$-dimensional manifold $N$, and the operator $A$ to be a pseudodifferential operator ($\Psi$DO). Furthermore, the noise in (\ref{contmodel1}) is modelled as $\delta\hspace{.2mm}\varepsilon(x)$, where $\delta>0$ is the noise amplitude and $\varepsilon=\varepsilon(x)$ is a realization of normalised Gaussian white noise $W(x)$. 

Rigorous treatment of white noise on $N$ is based on generalized functions (distributions). We denote the pairing of a distribution $f\in \mathcal D^\prime(N)$ and
 a test function $\phi\in C^\infty(N)$ by $\bra f,\phi\cet$.
  Let $(\Omega,\Sigma,\mathbb P)$ be a probability space.  
A random generalized function $V=V(x,\omega)$, where $x\in N$ and $\omega\in \Omega$, on $N$ 
is a measurable map $V:\Omega\to \mathcal D^\prime(N)$. Below, following the tradition used in study of stochastical processes, we often omit the $\omega$ variable and just
denote a random generalized function by $V(x)$.

White noise $W(x)$ is a random generalized function on $N$ such that the inner products
$\bra  W,\phi\cet $ are Gaussian random variables for all $\phi\in C^\infty(N)$,
 $\E W=0$, and
 \beq\label{gauss infinite}
\E\bigg(\bra W,\phi\cet \bra W,\psi\cet\bigg)=\bra \phi,C_W \psi\cet_{L^2(N)}\quad
\hbox{for }\phi,\psi\in C^\infty(N).
 \eeq
The covariance operator $C_W$ of  Gaussian white noise is the identity operator. Then $W$ can be considered
 as a function $W:\Omega\to \mathcal D^\prime(N)$ where $\Omega$ is the probability space. A realization of $W$ is the generalized function $x\mapsto W(x,\omega)$ on $N$ with a fixed $\omega\in \Omega$.

Below, we consider the case when 
\beq
\label{Pk definition}P_k(f)= (\bra f,\phi_j\cet)_{j=1}^k,
\eeq
where $\phi_j\in C^\infty (N)$ are such that 
$(\phi_j)_{j=1}^\infty$ is an orthogonal basis  in  $L^2(N)$. 
Then $P_k(\delta\hspace{.2mm}\varepsilon(x))=\delta\hspace{.2mm}\textbf{e}\in\R^k$ with $\textbf{e}$ as above.
For example, when $N$ is a $d$-dimensional torus, $P_k$ can be  the truncation of the Fourier series.
See  Section \ref{sec:Gamma convergence} for a   detailed discussion on  discrete and continuous noise models.


Now we can state the main motivation behind this study. 
The probability density
function of $W$ is often {\it formally} written in the form  
\begin{equation}\label{formaldensity}
\pi_W(w)=c\exp(-\|w\|_{L^2(N)}^2/2).
\end{equation} 
However, despite formula (\ref{formaldensity}), the realizations 
of the  white Gaussian noise are almost surely not in $L^2(N)$. 
Thus we cannot use formula (\ref{contTikh00}) when the error  in the measurement $m$ is white Gaussian noise.} Let us illustrate this ``white noise paradox'' by a simple example. 
\bigskip

\noindent
{\bf Example 1.}
{\em 
Let $W$ be {\newtext normalized Gaussian white noise defined on the $d$-dimensional torus 
$\T^d=(\R/(2\pi \Z))^d$.}
The Fourier coefficients of $W$ are normally distributed with variance one, that is, $\bra W,e_{\vec{\ell}}\cet\sim N(0,1)$, where $e_{\vec{\ell}}(x)=e^{i\vec{\ell}\cdot x}$ and $\vec{\ell}\in \mathbb{Z}^d$.
Hence
\ba
\E\|W\|_{L^2(\T^d)}^2=\sum_{\vec{\ell}\in \mathbb{Z}^d}\E|\bra W,e_{\vec{\ell}}\cet|^2=\sum_{\vec{\ell}\in \mathbb{Z}^d}1=\infty.
\ea
This implies that $W\in L^2(\T^d)$ with probability zero. 
However, when $s<-d/2$ 
\begin{equation}\label{sum_example1}
\E\|W\|_{H^{s}(\T^d)}^2=\sum_{\vec{k}\in \mathbb{Z}^d}(1+|\vec{\ell}|^2)^{s}\E|\bra W,e_\ell\cet|^2<\infty
\end{equation}
and hence $W$ takes values in $H^{s}(\T^d)$ 
{\newtext almost surely (that is, with probability one).}

On the other hand \cite[Theorem 2]{rozanov1971} implies that if $\|W\|_{H^{s}(\T^d)}^2<\infty$ almost surely then $\E\|W\|_{H^{s}(\T^d)}^2<\infty$ which yields $s<-d/2$. This concludes that the realisations of white noise $W$ are almost surely in the space $ H^{s}(\T^d)$
if and only if  $s<-d/2$. In particular for $s \geq -d/2$ the function $x\mapsto W(x,\omega)$ is in $H^{s}(\T^d)$ only when $\omega\in\Omega_0\subset\Omega$ where $\mathbb{P}(\Omega_0)=0$. }\\

Even though the previous example is proven in $\T^d$ we note that the same result is valid in all open bounded subsets $D\subset\R^d$. 

\subsection{Main result}
{\newtext Let us
again consider a general closed $d$-dimensional  Riemannian manifold $(N,g)$
and let $\Delta=\Delta_g$ be the Laplace operator on $N$.} Furthermore, let $A$ be a pseudodifferential operator.
Consider the following  measurement model:
\begin{align}\label{continuous}
m = Au + \delta\hspace{.2mm}\varepsilon,
\end{align}
where $\varepsilon\in H^s(N)$ with $s<-d/2$ is a realization of white noise. 

The pseudodifferential operator $A$ can be, for example,
\begin{equation*}
Au(x) = \int_N \mathcal{A}(x,z)u(z)dz
\end{equation*}  
where $\mathcal{A}\in C^\infty \big((N\times N)\backslash \text{diag}(N)\big)$ 
and in an open neighbourhood \linebreak
$U\opensubset N\times N$ of the diag$(N)=\{(x,x);\ x\in N\}$, we have 
\begin{equation*}
\mathcal{A}(x,z) = \frac{b(x,z)}{d_g(x,z)^p}, \quad (x,z)\in U
\end{equation*}
where $d_g$ is a distance function, $p<d$, $b\in C^\infty(U)$ and $b(x,x)\not=0$. In this case $\mathcal{A}$ is a pseudodifferential operator of order $-d+p<0$.

Let us now modify formula (\ref{contTikh00}) to arrive at something useful for white Gaussian noise. Expand the data fidelity term like this: $\|m-A u\|_{L^2(N)}^2=\|A u\|_{L^2(N)}^2- 2\langle m,A u \rangle +\|m\|_{L^2(N)}^2$. Simply omitting the ``constant term'' $\|m\|_{L^2(N)}^2$  leads to the definition
\begin{equation}\label{contTikh0}
  T_\alpha(m):=\argmin_{u\in H^r(N)}
  \big\{\|A u\|_{L^2(N)}^2 - 2\langle m,A u \rangle +
 \alpha\|u\|_{H^r(N)}^2 \big\},
\end{equation}
where we can interpret $\langle m, A u \rangle$ as a suitable duality pairing instead of $L^2(N)$ inner product. When {\newtext $A$ is a pseudodifferential operator of order $-t<s+r$, we can define $\langle m, Au \rangle = \langle m, Au \rangle_{H^{s}(N)\times H^{-s}(N)}$. 
}
  
It is well-known that the solution of the finite-dimensional problem (\ref{discrTikh1}) can be calculated using the following  formula:
\begin{equation}\label{discrTikh2}
  T_\alpha(\finm)=(\finA^T\!\finA+\alpha \finL^T\!\finL)^{-1}\finA^T \finm.
\end{equation}
The regularized solution of the continuous problem (\ref{contTikh0}) is
\begin{equation}\label{contTikh1}
  T_\alpha(m)=( A^*A+ \alpha(I-\Delta)^r)^{-1}A^*m.
\end{equation}
The regularization parameter is chosen to be a function of the noise amplitude: $\alpha(\delta)=\alpha_0\delta^\kappa,$ where $\alpha_0>0$ is a constant and $\kappa>0$.
We will now formulate the main theorem of this paper, concerning the continuous regularized solution (\ref{contTikh1}).  

\begin{theorem}\label{speedelliptic}
Let  $N$ be a $d$-dimensional closed manifold and $u\in H^r(N)$ with $r\geq 0$. Here  $\|u\|_{H^r(N)}:=\|(I-\Delta)^{r/2}u\|_{L^2(N)}$. Let $\varepsilon\in H^s(N)$ with some $s<-d/2$ and consider the measurement 
 \beq\label{eq: measurement}
 m_\delta=Au+\delta\hspace{.2mm}\varepsilon,
 \eeq
where $A\in\Psi^{-t}$, is an elliptic pseudodifferential operator of order $-t$ {\newtext on the manifold $N$ with $t>\max\{0, -s-r\}$
 and $\delta\in \R_+$}. Assume that $A:L^2(N)\to L^2(N)$ is injective. The regularization parameter is chosen to be $\alpha(\delta)=\alpha_0\delta^\kappa,$ where $\alpha_0>0$ is a constant and $\kappa>0 $.

Take $s_1 \leq s-t+2(t+r)/\kappa$. Then the following convergence takes place in $H^{s_1}(N)$ norm: 
\ba
\lim_{\delta\to 0} T_{\alpha(\delta)}(m_\delta)=u.
\ea  
Furthermore, we have the following estimates for the speed of convergence:
\begin{itemize}
\item[(i)] If $s_1\leq s-t$ then 
\begin{align*}
\|T_{\alpha(\delta)}(m_\delta)-u\|_{H^{s_1}}\leq C\max\{\delta^{\frac{\kappa(r-\zeta)}{2(t+r)}},\delta\}.
\end{align*}  
\item[(ii)] If $s-t\leq s_1 < s-t+2(t+r)/\kappa$ then 
\begin{align*}
\|T_{\alpha(\delta)}(m_\delta)-u\|_{H^{s_1}}\leq C\max\{\delta^{\frac{\kappa(r-\zeta)}{2(t+r)}},\delta^{1+\frac{\kappa(s-t-s_1)}{2(t+r)}}\}.
\end{align*}
\end{itemize} 
Above we have $\zeta=\max\{s_1,-r-2t\}$. 
\end{theorem}
Notice that in case (i) $\omega\geq 0$ and in case (ii) $\omega \leq 0$.
The different convergence speeds (i) and (ii) show the trade-off between smoothness of the space and the speed of convergence. In case (i) we get better convergence rates but in case (ii) we can use a stronger norm.  In section \ref{sec:deblurring} we give two counterexamples  to  show that even though $u\in H^r$ and $T_{\alpha(\delta)}(m_\delta)\in H^r$ the regularized solution does not converge to the real solution in $H^r$ norm.

\subsection{Literature review}
There are two main ways in inverse problems literature for modelling noise. The first approach based on the deterministic regularization techniques 
is to assume that the noise is deterministic and small. In that case one has a norm estimate of the noise and can study what happens when $\|
{\rm noise}\|_{L^2}\to 0$. This approach was originated by Tikhonov \cite{Tikhonov1,Tikhonov2}, and studied in depth in \cite{Burger,EHN,Groetsch,Kirch,Philips,Morozov, Tikhonov3}.  
%
The second  approach to handling the noise is based on statistical point of view. The statistical modeling 
of noise in the inverse problems started in the early papers of \cite{Franklin,fitzpatrick,Sudakov, 
tarantolabook} and it is notable
that with this approach one needs not assume smallness of the noise. For some recent references of the frequentist view of statistical problems see \cite{hohage1,hohage2,Mair,MatPer}. Another statistical way to study
inverse problems with random noise is based on Bayesian approach where $m,\ x$ and $\varepsilon$ are considered to be realizations of random variables, see \cite{EkiDaniela2,hanson87,Kaipio,Lasanen1,Lasanen2,Lasanen3,LSS,LS,LPS,Pikkarainen1,
Pikkarainen2,Pi}.

The  deterministic regularization and statistical approaches differ both in assumptions and techniques.
This paper aims to bridge the gap between them. Our results are closely
related to earlier studies of Eggermont, \mbox{LaRiccia}, and  Nashed \cite{nashed1,nashed2,nashed3},
who studied weakly bounded noise.
They assume that the noise is a $L^2$-function and discuss regularization techniques when the noise tends to zero in the weak topology of $L^2$. This kind of relaxed assumption of noise covers small low frequency noise and large high frequency noise.
However, even though  $\delta\hspace{.2mm}\varepsilon$ tends to zero in weak sense as $\delta\to 0$ when $\eps$ is a realization
of the normalized white noise, this type of noise lies outside the definition of the weakly bounded noise as 
$\eps$ is not almost surely $L^2$-valued.

A related approach of smoothing the noise before the analysis is described in \cite{Mathe1,Mathe2}.
A similar regularization method where no smoothness of the operator $A$ is assumed, but instead the regularization method is modified, is studied in \cite{Egger}. Another possible approach to deal with white noise is to first perform a data projection step and then proceed to Tikhonov regularization \cite{KR,KMR}. Also, Hohage and Werner have earlier studied inverse problems taking
into account the fact that white noise is not square-integrable in \cite{hohage3}.

Our new results are different from all of those previous studies. Our approach aims to study the effect of the continuous white noise not being an $L^2$ function in Tikhonov regularisation (\ref{contTikh00}) instead of modifying the problem by altering the regularisation method or assumptions.

%
%

\section{Analysis of the translation-invariant case}\label{sec:translationexample}  

Before giving the general proof of Theorem \ref{speedelliptic} in Section \ref{sec:generalproof} we  motivate the proof by proving a similar kind of lemma for translation-invariant case. 

The regularized solution we are studying is of the form
\begin{equation*}
  T_{\alpha(\delta)}(m):=\argmin_{u\in H^r(\T^d)}
  \big\{\|A u\|_{L^2(\T^d)}^2 - 2\langle m,A u \rangle +\alpha\|(I-\Delta)^{r/2}u\|_{L^2(\T^d)}^2 \big\},
\end{equation*}
where $\alpha(\delta)=\alpha_0\delta^\kappa$, for some constant $\alpha_0$ and $\kappa>1$. As mentioned before solution to this is 
\begin{equation}\label{regsolutio}
  T_{\alpha(\delta)}(m)=( A^*A+ \alpha(I-\Delta)^r)^{-1}A^*m.
\end{equation}

Let us consider the case when $\alpha=\alpha_0\delta^2$, where $\alpha_0$ is a constant, and $A$ is an elliptic pseudodifferential 
operator of order $-t<0$ that commutes
with translations. Then, in $L^2(\T^d)$ we have that  $B=A^*A\geq c_1(I-\Delta)^{-t}$.
As $A$ and $B$ commute with translations they are Fourier multipliers,
\ba
\widehat{Au}(n)=a(n)\widehat{u}(n)
\ea
and since $A$ is elliptic there is $n_0>0$ so that
\ba
c_1|n|^{-t}\leq  |a(n)|\leq c_2|n|^{-t},\quad\hbox{for }|n|>n_0.
\ea
The symbol $z_\delta$ of $Z_\delta= A^*A+\alpha_0\delta^2(I-\Delta)^r$ 
is
\ba
z_\delta(n)= |a(n)|^2+{\alpha_0\delta^2}(1+n^2)^r
\ea
and thus
\ba
z_\delta(n)\geq \max( |a(n)|^2, \alpha_0\delta^2(1+n^2)^r).
\ea
If $0<\beta<\frac 12$ and $|n|>n_0$
\ba
z_\delta(n)&\geq& |a(n)|^{2(1-\beta)}(\alpha_0 \delta^2(1+n^2)^r)^\beta  
\\
&\geq& c_3|n|^{-2(1-\beta)t+2r\beta }\delta^{2\beta} \alpha_0^{\beta}.
\ea
Now when $s<-d/2$ we have
\ba
\varepsilon\in H^s(\T^d).
\ea
Thus writing   
\begin{equation}\label{u2vw}
\begin{split}
T_{\alpha(\delta)}(m_\delta) = & \bigg( A^*A+ {\alpha_0\delta^2}(I-\Delta)^r\bigg)^{-1}(A^*(Au+\delta\hspace{.2mm}\varepsilon))\\
= &\bigg( A^*A+ {\alpha_0\delta^2}(I-\Delta)^r\bigg)^{-1}A^*Au+\\
& +\bigg( A^*A+ {\alpha_0\delta^2}(I-\Delta)^r\bigg)^{-1}A^*(\delta\hspace{.2mm}\varepsilon)
\end{split}
\end{equation}
we see that
\ba 
 T_{\alpha(\delta)}(m_\delta)  &=&v_\delta+ w_\delta
\ea
where
\ba
& &\widehat v_\delta (n)=\frac 1{z_\delta(n)}|a(n)|^2\widehat {u}(n),
\\
& &\widehat w_\delta (n)=\frac 1{z_\delta(n)}\overline{a(n)} \widehat\varepsilon(n)\delta .
  \ea
  Here,
\ba
\bigg|\frac 1{z_\delta(n)}|a(n)|^2\bigg|\leq 1,\quad 
\lim_{\delta\to 0}\frac 1{z_\delta(n)}|a(n)|^2=1
\ea
and thus if $u\in H^r(\T^d)$ by dominated convergence theorem
\ba
\lim_{\delta\to 0}v_\delta=u,\quad \hbox{in }H^r(\T^d).
\ea
Above the limit speed of convergence can be analysed using the
standard regularization theory \cite{EHN} and the fact that 
\ba
 v_\delta  &=&\bigg( A^*A+ {\alpha_0\delta^2}(I-\Delta)^r\bigg)^{-1}A^*Au\\
  &=&u-
{\alpha_0\delta^2} \bigg( A^*A+ {\alpha_0\delta^2}(I-\Delta)^r\bigg)^{-1}(I-\Delta)^ru. 
   \ea
We can use the fact that 
$Z_\delta=A^*A+ {\alpha_0\delta^2}(I-\Delta)^r\geq  {\alpha_0\delta^2}(I-\Delta)^r$ and write 
\ba
\|Z_\delta^{-1/2} (I-\Delta)^ru\|_{L^2}&\leq&
(\alpha_0\delta^2)^{-1/2}\|  (I-\Delta)^{r/2}u\|_{L^2}\\
&\leq&
(\alpha_0\delta^2)^{-1/2}\| u\|_{H^r}.
\ea
We also have the inequality $Z_\delta\geq A^*A\geq c_1(I-\Delta)^{-t}$.
When $r>0$ we can define $\eta=t/(2r+2t)$ and $\gamma=r/(2r+2t)$ 
so that, $\gamma+\eta=1/2$, $t\gamma-r\eta=0$. We get

\begin{align*}
\| Z_\delta^{-1} &(I-\Delta)^ru\|_{L^2}
=\| Z_\delta^{-\gamma-\eta-1/2} (I-\Delta)^ru\|_{L^2}\\
&\leq
(\alpha_0\delta^2)^{-1/2}\| Z_\delta^{-\gamma-\eta} (I-\Delta)^{r/2}u\|_{L^2}
\\
&\leq
(\alpha_0\delta^2)^{-1/2}\| (c_1(I-\Delta)^{-t})^{-\gamma}  (\alpha_0\delta^2 (I-\Delta)^r)^{-\eta} (I-\Delta)^{r/2}u\|_{L^2}
\\
&\leq
 c_1^{-\gamma}  (\alpha_0\delta^2)^{-\eta-1/2} \| (I-\Delta)^{r/2}u\|_{L^2}
\\
&\leq 
 c_1^{-\gamma}  (\alpha_0\delta^2)^{-\eta-1/2} \| u\|_{H^r}.
\end{align*}
Hence we obtain

\ba
\begin{split}
\|{\alpha_0\delta^2} Z_\delta^{-1}(I-\Delta)^ru\|_{L^2} & \leq c_1^{-\gamma}  (\alpha_0\delta^2)^{1/2-\eta} \| u\|_{H^r}\\
& = c_1^{-\gamma}  \delta^{\frac{r}{t+r}} \| u\|_{H^r}& .
\end{split}
\ea
On the other hand,

\ba
\bigg|\frac 1{z_\delta(n)}\overline{a(n)} \delta\bigg| &\leq&
\frac 1{c_3|n|^{-2(1-\beta)t+2r\beta}\delta^{2\beta} \alpha_0^{\beta}} c_2|n|^{-t}\delta\\
&\leq&
c_4|n|^{(1-2\beta)t-2r\beta}\delta^{1-2\beta} \alpha_0^{-\beta}.
\ea
Hence
\ba
\|w_\delta\|_{H^{s_1}(\T^d)}\leq c_5\delta^{1-2\beta}  
\ea
where $s_1\leq s-(1-2\beta)t+2r\b$. Because we proved the convergence of $v_\delta$ in $L^2$  we have to have $s_1\leq 0$. This is true at least when $s\leq -r$.
Thus adding the above results together we can formulate the next lemma.

\begin{lemma} Let $u\in H^r(\T^d)$, $r>0$, be Gaussian distributed, $\varepsilon\in H^s(\T^d)$, $s<\max\{-d/2,-r\}$, and
 \ba
 m_\delta=Au+\delta\hspace{.2mm}\varepsilon
 \ea
where $A:H^r(\T^d)\to H^{r+t}(\T^d)$, $t>\max\{0, -s-r\}$, is an elliptic pseudodifferential operator of order $-t<0$ that commutes with translations. We assume that $\alpha(\delta)=\alpha_0\delta^2$. Then for the regularized solution  
\begin{equation*}
T_{\alpha(\delta)}(m)=( A^*A+ \alpha(I-\Delta)^r)^{-1}A^*m
\end{equation*}
of  $u$ we have

\ba
\lim_{\delta\to 0} T_{\alpha(\delta)}(m_\delta)=u,\quad \hbox{in }H^{s_1}(\T^d)
\ea
where $s_1\leq s-(1-2\beta)t+2r\b\leq 0$ and $0<\beta<1/2$. Furthermore we have the following estimate of the speed of convergence 
\ba 
\|T_{\alpha(\delta)}(m_\delta)-u\|_{H^{s_1}}\leq C\max\{\delta^{\frac{r}{t+r}},\delta^{1-2\beta}\}.
\ea
\end{lemma}

\noindent\textbf{Proof.} The convergence is immediate consequence of the above results. For the convergence speed we get 
\ba 
\|T_{\alpha(\delta)}(m_\delta)-u\|_{H^{s_1}} &=& \|-\alpha_0\delta^2 Z_\delta^{-1}(I-\Delta)^ru+w_\delta\|_{H^{s_1}}\\
&\leq & \|\alpha_0\delta^2 Z_\delta^{-1}(I-\Delta)^ru\|_{L^2}+\|w_\delta\|_{H^{s_1}}\\
&\leq & C_1(\alpha_0\delta^2)^{1/2-\eta}+C_2\delta^{1-2\b}\alpha_0^{-\beta}\\
&\leq & C_3\max\{\delta^{\frac{r}{t+r}},\delta^{1-2\beta}\}
\ea
where $Z_\delta=A^*A+\alpha_0\delta^2(I-\Delta)^r$ and $\eta=t/2(t+r)$.\hfill  $\square$

\section{Proof of the main theorem}   \label{sec:generalproof}

Here we study the general case where $A$ is an elliptic pseudodifferential operator of order $-t<0$. We denote $H^s(N)=H^s$ and $L^2(N)=L^2$ where $N$ is a closed manifold and $\dim N=d$. As in the previous example we have
\begin{equation} \label{uparts}
\begin{split}
T_{\alpha(\delta)}(m_\delta) & = Z_\delta^{-1}A^*Au+Z_\delta^{-1}A^*(\delta\hspace{.2mm}\varepsilon)\\
& = u-\alpha Z_\delta^{-1}(I-\Delta)^ru+Z_\delta^{-1}A^*(\delta\hspace{.2mm}\varepsilon)
\end{split}
\end{equation}
where $Z_\delta=A^*A+\a(I-\Delta)^r$. 

First we will show that $B=A^*A$ is invertible. 
We define $A^*:L^2(N)\to L^2(N)$ as the adjoint of an operator $A:L^2(N)\to L^2(N)$. 
We assume that $A:L^2(N)\to L^2(N)$ is one-to-one. If $Bu=0$ then 
\ba 
0=\langle A^*Au,u \rangle_{L_2}=\langle Au,Au \rangle_{L_2}=\|Au\|_{L_2}^2
\ea
which implies $Au=0$ and furthermore $u=0$. Thus the operator $B:L^2(N)\to H^{2t}(N)$ is one-to-one. 

Next we recall the fact that an elliptic operator $B\in \Psi^{-2t}(N)$ is a Fredholm operator and $\ind(B)=0$ (\cite{hormander3} Theorem 19.2.1). Indeed index of a Fredholm operator $B$ is 
\begin{align}\label{index}
\ind (B)=\dim(\Ker B)-\dim(\coker B).
\end{align} 
If $K:L^2(N)\to H^{2t}(N)$ is compact and $B^{adj}:(H^{2t})^*(N)=H^{-2t}(N)\to (L^2)^*(N)=L^2(N)$ is the adjoint of the operator $B:L^2(N)\to H^{2t}(N)$ then 
\ba 
\ind (B+K)=\ind(B)=-\ind(B^{adj})
\ea
Define $B_s:H^s(N)\to H^{s+2t}(N)$ as an extension of $B:C^\infty(N)\to C^\infty(N)$ and show that $\ind (B_s)=\ind(B_0)$ for all $s$. Define
\ba 
P=(I-\Delta_g)^{s/2}B_s(I-\Delta_g)^{-s/2}:L^2(N)\to H^{2t}(N).
\ea
We can write $P=B_0+K_1$ where $K_1:L^2(N)\to H^{2t}(N)$ is compact. 
Now
\ba 
\ind (B_0) &=& \ind (P)\\
&=& \ind (I-\Delta_g)^{s/2}+\ind (B_s) +\ind (I-\Delta_g)^{-s/2} \\
&=& \ind (B_s).
\ea
Because $B-B^{adj}:L^2(N)\to H^{2t}(N)$ is compact
we can write 
\ba 
\ind (B:L^2(N)\to H^{2t}(N)) &=& -\ind (B^{adj}:H^{-2t}(N)\to L^2(N))\\
&=& -\ind (B^{adj}:L^2(N)\to H^{2t}(N))\\
&=& -\ind (B:L^2(N)\to H^{2t}(N))
\ea
and hence we see that $\ind (B:L^2(N)\to H^{2t}(N))=0$. Using this, the knowledge that $B$ is one-to-one and (\ref{index}) we get
\ba 
0=\dim(\Ker B)=\dim(\coker B)
\ea 
which means that $B$ is also onto. Thus we have shown that there exist $B^{-1}:H^{2t}(N)\to L^2(N)$.

Next we will examine $\Psi$DOs that depend on spectral variable $\lambda=(\alpha_0\delta^\kappa)^{-1}$. For the general theory see \cite{shubin}.     
The symbol class $S^m_p(\R^d\times\R^d,\R_+)$ consist of the functions $a(x,\xi,\lambda)$ such that  
\begin{enumerate}
  \item $a(x,\xi,\lambda_0)\in C^\infty(\R^d\times\R^d)$ for every fixed $\lambda_0\geq0$ and
  \item for arbitrary multi-indices $\alpha$ and $\beta$ and for any compact set $K\subset \R^d$ there exist constants $C_{\alpha,\beta,K}$ such that 
  \[|\partial^\alpha_\xi\partial^\beta_x a(x,\xi,\lambda)|\leq C_{\alpha,\beta,K}(1+|\xi|+|\lambda|^{1/p})^{m-|\alpha|}\]  
  for $x\in K$, $\xi\in \R^d$ and $\lambda\geq0$.
\end{enumerate}
We consider the pseudodifferential operators $A_\lambda:\mathcal D^\prime (N)\to \mathcal D^\prime(N)$ depending on the parameter $\lambda$. To define such operators, one considers local coordinates $Y:U\to \R^d$ of the manifold $N$, where
we emphasize that the set $U\subset N$ does not need to be connected (see \cite[Sect.\ I.4.3]{shubin}).
A bounded linear operator  $A_\lambda:\mathcal D^\prime (N)\to \mathcal D^\prime(N)$,
depending on the parameter $\lambda$, is a  pseudodifferential operator with spectral variable $\lambda$ if 
 for any local coordinates $Y:U\to \R^d$ of manifold $N$, $U\subset N$, there
 is a symbol $a\in S^m_p(\R^d\times\R^d,\R_+)$ such that
for  $u\in C^\infty_0(U)$  we have 
\begin{align*}
(A_\lambda u)(Y^{-1}(x))=\int_{V\times \R^d} e^{i(x-y)\cdot\xi} a(x,\xi,\lambda)u(Y^{-1}(y))dyd\xi,\quad x\in V,
\end{align*}
where $V=Y(U)\subset \R^d$.
In this case we will write 
\[A_\lambda\in \Psi_p^m(N,\R_+),\]
and say that in local coordinates $Y:U\to \R^d$ the operator $A$ has the symbol $a(x,\xi,\lambda)\in S^m_p(\R^d\times\R^d,\R_+)$. 
If for all compact sets $K\subset \R^d$ there are constants $C_1,C_2,R>0$  such that the symbol $a(x,\xi,\lambda)\in S^m_p(\R^d\times \R^d,\R_+)$ satisfies 
\ba 
C_1(|\xi|+|\lambda|^{1/p})^m\leq |a(x,\xi,\lambda)|\leq C_2(|\xi|+|\lambda|^{1/p})^m,
\ea 
for $|\xi|+|\lambda|\geq R$ and $x\in K$, we say that $a$ is hypoelliptic with parameter and denote $a(x,\xi,\lambda)\in HS^m_p(\R^d\times \R^d,\R_+)$. We will denote by $H\Psi^m_p(N,\R_+)$ the class of $\Psi$DOs depending on the parameter $\lambda$ whose symbol  in all local coordinates belongs in $HS^m_p(\R^d\times \R^d,\R_+).$

We want to prove that 
\ba 
F_\lambda=( A^*A)^{-1}(I-\Delta)^r+\lambda I
\ea    
is invertible. Operator $F_\lambda\in \Psi^{2(t+r)}(N)$ is elliptic since $( A^*A)^{-1}(I-\Delta)^r\in \Psi^{2(t+r)}(N)$ is elliptic and $\lambda\ I\in \Psi^0(N)$. Denote $Q=( A^*A)^{-1}(I-\Delta)^r$ and its symbol $q(x,\xi)\in S^{2(t+r)}(N)$. Then for the symbol $\sigma(F_\lambda)(x,\xi)=q(x,\xi)+\lambda$ of the operator $F_\lambda$ we have in 
compact subsets $K$ of  any local coordinates
\ba 
|\partial^\a_\xi \partial^\b_x (q(x,\xi)+\lambda)| \leq C_{\a,\b,K}(1+|\xi|+|\lambda|^{1/(2(t+r))})^{2(t+r)-|\a|},\ x\in K.
\ea
By (\cite{shubin} Theorem 9.2.) there exist $R>0$ such that for $|\lambda|\geq R$ the operator $F_\lambda\in H\Psi^{2(t+r)}_{2(t+r)}(N,\R_+)$ is invertible with
\ba 
F^{-1}_\lambda\in H\Psi^{-2(t+r)}_{2(t+r)}(N,[R,\infty)).
\ea

Now we have shown that the operator $Z_\delta^{-1}$ can be rewritten  
\begin{equation} 
Z_\delta^{-1} = \lambda\bigg(( A^*A)^{-1}(I-\Delta)^r+\lambda\bigg)^{-1}( A^*A)^{-1}
\end{equation}
where $\lambda=(\alpha_0\delta^\kappa)^{-1}$.

We denote by $\| F_\lambda\|_{s,s-\ell}$ the norm of $F_\lambda:H^s(N)\to H^{s-\ell}(N)$ where $s,\ell\in \R$. We have the following norm estimates for  $F_\lambda\in \Psi^{m}_p(N,\R_+)$ when $ \ell\geq m$ and $\lambda$ large enough  
\begin{align}
\| F_\lambda\|_{s,s-\ell}\leq C_{s,l}(1+|\lambda|^{1/p})^m, \quad\quad \text{if}\quad \ell\geq 0 \label{Shubin1}\\
\| F_\lambda\|_{s,s-\ell}\leq C_{s,l}(1+|\lambda|^{1/p})^{-(\ell-m)}, \quad\quad \text{if}\quad \ell\leq 0. \label{Shubin2}
\end{align}

We can rewrite (\ref{uparts}) 
\begin{align}\label{uparts2}
T_{\alpha(\delta)}(m_\delta) & = u-\alpha Z_\delta^{-1}(I-\Delta)^ru+\alpha_0^{-1}\delta^{1-\kappa} F_\lambda^{-1}(A^*A)^{-1}A^*\varepsilon.
\end{align}

For the third term on the right hand side of (\ref{uparts2}) we have  $(A^*A)^{-1}A^*:H^s(N)\to H^{\tilde{s}}(N)$, $\tilde{s}=s-t<-d/2$ and hence we get 
\ba 
\|F_\lambda^{-1}(A^*A)^{-1}A^*\varepsilon\|_{H^{s_1}} &=& \|F_\lambda^{-1}\|_{\tilde{s},s_1}\|(A^*A)^{-1}A^*\varepsilon\|_{H^{\tilde{s}}}\\
&\leq& C\|F^{-1}_\lambda\|_{\tilde{s},s_1}.
\ea
Above $F_\lambda^{-1}\in \Psi^m_p(N,\R_+)$, where $m=-2(t+r)$ and $p=2(t+r)$. 

First we study the case when $s_1\leq\tilde{s}=s-t<-d/2$. Inequality (\ref{Shubin1}) gives us the norm estimate 
\ba 
\|F_\lambda^{-1}\|_{\tilde{s},\tilde{s}-\ell}\leq C(1+|\lambda|^{1/p})^m=C(1+\delta^{-\kappa/p})^m,
\ea 
 where $\ell=\tilde{s}-s_1\geq 0$. Now clearly $\ell\geq m=-2(t+r)$. Because we want $\delta^{1-\kappa}\|F_\lambda^{-1}(A^*A)^{-1}A^*\varepsilon\|_{H^{s_1}}\to 0$ when $\delta\to 0$ we have to require that 
\ba
-\frac{\kappa m}{p}=\kappa>\kappa-1.
\ea
which is true for all $\kappa,t>0$ and $r,\ell\geq 0$. 
 
 When $\tilde{s}\leq s_1$ we can use (\ref{Shubin2})
\ba
\|F_\lambda^{-1}\|_{\tilde{s},\tilde{s}-\ell}\leq C(1+|\lambda|^{1/p})^{-(\ell-m)}=C(1+\delta^{-\kappa/p})^{-(\ell-m)}, 
\ea
where $0 \geq \tilde{s}-s_1=\ell \geq m=-2(t+r)$ if $s_1\leq s+t+2r$. For convergence we need 
\ba 
\frac{\kappa(\ell-m)}{p}=\kappa\bigg(1+\frac{\ell}{p}\bigg)>\kappa-1
\ea
that is $0\geq\ell > -2(t+r)/\kappa$. For $s_1$ we get $s-t\leq s_1 \leq s-t+2(t+r)/\kappa$.

When $r>0$ the convergence of $\alpha Z_\delta^{-1}(I-\Delta)^ru$ could be shown the same way as in the previous example. Next we will show the convergence also in the case $r=0$ and improve the convergence rate by proving the convergence in $H^\zeta$ instead of $L^2$.

Assume that $r \geq 0$ and denote $\zeta =-r-\theta\geq s_1$. We need to find such $\eta\geq0$ and $\gamma\geq0$ that  $\gamma+\eta=1$ and $t\gamma-r\eta-\theta/2=0$. Define $\eta=(2t-\theta)/2(t+r)$ and $\gamma=(2r+\theta)/2(t+r)$, where $\theta\leq2t$. Using the inequalities $Z_\delta=A^*A+ {\alpha}(I-\Delta)^r\geq  {\alpha}(I-\Delta)^r$ and $T\geq A^*A\geq c_1(I-\Delta)^{-t}$ we get
\ba
\begin{split}
\|{\alpha} Z_\delta^{-1}(I-\Delta)^ru\|_{H^{\zeta}} 
& \leq \alpha\|(c_1(I-\Delta)^{-t})^{-\gamma}(\alpha(I-\Delta)^{r})^{-\eta} (I-\Delta)^{r-\frac{r}{2}-\frac{\theta}{2}} u\|_{L^2}\\
& \leq c_1^{-\gamma}  \alpha^{1-\eta}\|(I-\Delta)^{t\gamma-r\eta -\frac{\theta}{2}}(I-\Delta)^{\frac{r}{2}}  u\|_{L^2}\\
& = c_1^{-\gamma} \delta^{\frac{\kappa(r-\zeta)}{2(t+r)}} \| u\|_{H^r}
\end{split}
\ea
where $\zeta =\max\{s_1,-r-2t\}$.

Adding the above results together we can prove Theorem \ref{speedelliptic}.
\medskip

\noindent\textbf{Proof.} Proof of Theorem \ref{speedelliptic}. The convergence is immediate consequence of the above results. Now when $s_1\leq\tilde{s}$ we have 
\ba   
\|T_{\alpha(\delta)}(m_\delta)-u\|_{H^{s_1}} 
&\leq &\a\| Z_\delta^{-1}(I-\Delta)u\|_{H^\zeta}+\alpha_0^{-1}\delta^{1-\kappa}\| F_\lambda^{-1}(A^*A)^{-1}A^*\varepsilon\|_{H^{s_1}}\\
&\leq & C_1\delta^{\frac{\kappa(r-\zeta)}{2(t+r)}}+C_2\delta^{1-\kappa+\frac{\kappa m}{p})}\\
&\leq & C_3\max\{\delta^{\frac{\kappa(r-\zeta)}{2(t+r)}},\delta\}.
\ea

If $\tilde{s}\leq s_1\leq s-t+2(t+r)/\kappa$ we get  
\ba  
\|T_{\alpha(\delta)}(m_\delta)-u\|_{H^{s_1}} 
&\leq & \a\| Z_\delta^{-1}(I-\Delta)^r u\|_{H^\zeta}+\alpha_0^{-1}\delta^{1-\kappa}\| F_\lambda^{-1}(A^*A)^{-1}A^*\varepsilon\|_{H^{s_1}}\\
&\leq & C_1\delta^{\frac{\kappa(r-\zeta)}{2(t+r)}}+C_2\delta^{1-\kappa+\kappa(1+\frac{\ell}{p})}\\
&\leq & C_3\max\{\delta^{\frac{\kappa(r-\zeta)}{2(t+r)}},\delta^{1+\frac{\kappa(s-t-s_1)}{2(t+r)}}\}.  
\ea
Above $\zeta=\max\{s_1,-r-2t\}$.

%
\hfill  $\square$

\section{A model problem: one-dimensional deblurring}\label{sec:deblurring}  

We consider a simple inverse problem to give flavour of results for the reader. Let $\T^2$ be the two-dimensional torus constructed by identifying parallel sides of the square $D=(0,1)^2\subset \R^2$; we model periodic
images as elements of function spaces over $\torus^2$. The {\mntext continuum model} is $m=Au+\noise$ with convolution operator $A$
defined by
\begin{equation}\label{A_conv}
  Au(x)=\int_{\T^2} \Phi(x-y)u(y)\,dy,
\end{equation}
{\newtext where $\Phi\in C(\T^2)$ is a point spread function that is given by the Schwartz kernel
of an elliptic pseudodifferential operator of the order $-t<-2$}.


\subsection{Divergence in $H^1$ norm}

Let us return to the translation-invariant case where we wrote the regularised solution in the form 
\ba 
 T_{\alpha(\delta)}(m_\delta)  &=&v_\delta+ w_\delta
\ea
where 
\begin{equation*}
v_\delta =  \bigg( A^*A+ {\alpha_0\delta^2}(I-\Delta)\bigg)^{-1}A^*Au
\end{equation*}
and
\ba 
w_\delta = \bigg( A^*A+ {\alpha_0\delta^2}(I-\Delta)\bigg)^{-1}A^*(\delta\hspace{.2mm}\varepsilon).
\ea
Now the Fourier transform of $w_\delta$ is 
\ba 
\widehat{w}_\delta(n) = \frac 1{z_\delta(n)}\overline{a(n)} \widehat\varepsilon(n)\delta.
\ea
Denote $I(\delta)=\{n\ |\ c_0\delta^2(1+n^2)\leq |a(n)|^2\leq c_1\delta^2(1+n^2)\}$. We get 
\ba 
\|w_\delta\|^2_{H^1} &\geq & \sum_{I(\delta)} (1+n^2)\bigg|\frac{\overline{a(n)}\delta}{|a(n)|^2+\alpha_0\delta^2(1+n^2)} \widehat\varepsilon(n)\bigg|^2\\
&\geq & \sum_{I(\delta)} \frac{|a(n)|^2}{|a(n)|^2+\alpha_0\delta^2(1+n^2)}\cdot\frac{(1+n^2)\delta^2}{|a(n)|^2+\alpha_0\delta^2(1+n^2)} |\widehat\varepsilon(n)|^2\\
&\geq & \sum_{I(\delta)} \frac{1}{1+\frac{\a_0}{c_0}}\cdot\frac{1}{c_1+\a_0} |\widehat\varepsilon(n)|^2.
\ea
We can chose $c_0,c_1$ so that $I(\delta)\not=\emptyset$ for all $0<\delta<\delta_0$. Now there exist $n(\delta)\in I(\delta)$ and $n(\delta)$ goes trough all $\{n\in \N\ |\ n\geq n_0\}$ when $\delta\to 0$. We see that 
\ba 
\limsup_{\delta\to 0}\|w_\delta\|^2_{H^1}\geq \limsup _{\delta\to 0} c_2 |\widehat\varepsilon(n(\delta))|^2 \geq c_2
\ea
almost surely since $\varepsilon$ is white noise.  Thus the solution $T_{\alpha(\delta)}(m_\delta)$ does not converge in $H^1$.
 
\subsection{Computational results}

Since the operator $A$ does not have a continuous inverse operator $L^2\to L^2$,  the condition number of the matrix approximation $\finA$ of the operator $A$ grows when the discretization is refined, i.e., when $n\ra\infty$ or $k\ra \infty$.
This is the very reason why regularization is need in the (numerical) solutions of the inverse problems

Next we demonstrate the above results numerically and 
 consider one-dimensional deblurring problem {\newtext  on the torus $\T^1=\R/\Z$,} 
\begin{equation*}
m = Au+\delta\hspace{.2mm}\varepsilon,
\end{equation*}
where $u\in H^1(\T^1)$ is the following piecewise linear function:
\begin{displaymath}
u = \left\{ \begin{array}{ll}
0 &  \textrm{when $0<x<0.3$ or $0.7<x<1$}\\
10x-3 & \textrm{when $0.3<x<0.4$}\\
1 & \textrm{when $0.4\leq x \leq 0.6$}\\
-10x+7 & \textrm{when $0.6<x<0.7$,}
\end{array} \right.
\end{displaymath}
$\varepsilon\in H^s$, $s<-1/2$ is white noise and $A$ is a $2$ times smoothing operator
\[(Au)(x)=\mathcal{F}^{-1}\big((1+|n|^2)^{-1}(\mathcal{F}u)(n)\big)(x).\]
Now solving $u$ from $Au(x)=m(x)$ corresponds to the solution of ordinary differential equation $(1-\partial_x^2)m(x)=u(x)$ so $A$ can be thought e.g. as a blurring operator. 

%

We assume $\alpha=\delta^{5/2}$ and thus the regularized solution is 
\ba 
u_\delta=T_{\alpha(\delta)}(m_\delta)=( A^*A+ \delta^{5/2}(I-\Delta)^r)^{-1}A^*m.
\ea

\begin{figure}[h]
\centering
\begin{subfigure}{.5\textwidth}
  \centering
  \includegraphics[height=5.2cm]{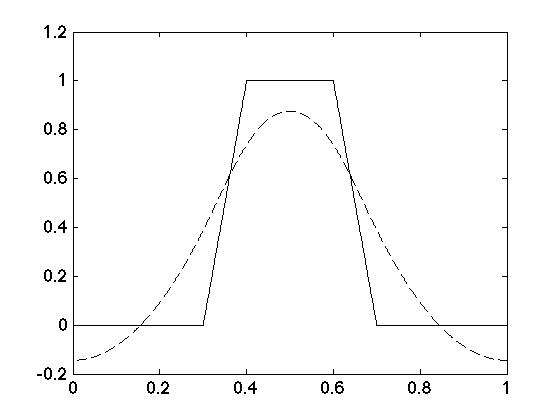}
\end{subfigure}%
\begin{subfigure}{.5\textwidth}
  \centering
  \includegraphics[height=5.2cm]{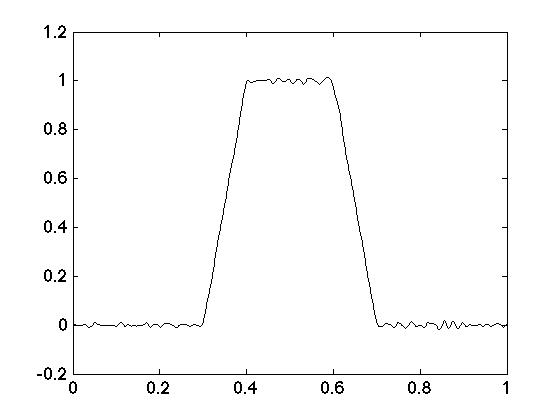}
\end{subfigure}
\caption{On the left the original piecewise linear function $u$ (solid line) and the noiseless data $m=Au$ (dashed line). On the right regularised solution $u_\delta$ when $\delta = 3,5*10^{-5}$.}
\label{fig: hat}
\end{figure}

Now Theorem \ref{speedelliptic} gives us 
\ba 
\lim_{\delta \to 0} \| u-u_\delta \|_{H^{s_1}} = 0
\ea
when $s_1<s-t+(t+r)/\kappa<-13/10$.

\begin{figure}[h!]
\begin{center}
\includegraphics[height=6.6cm]{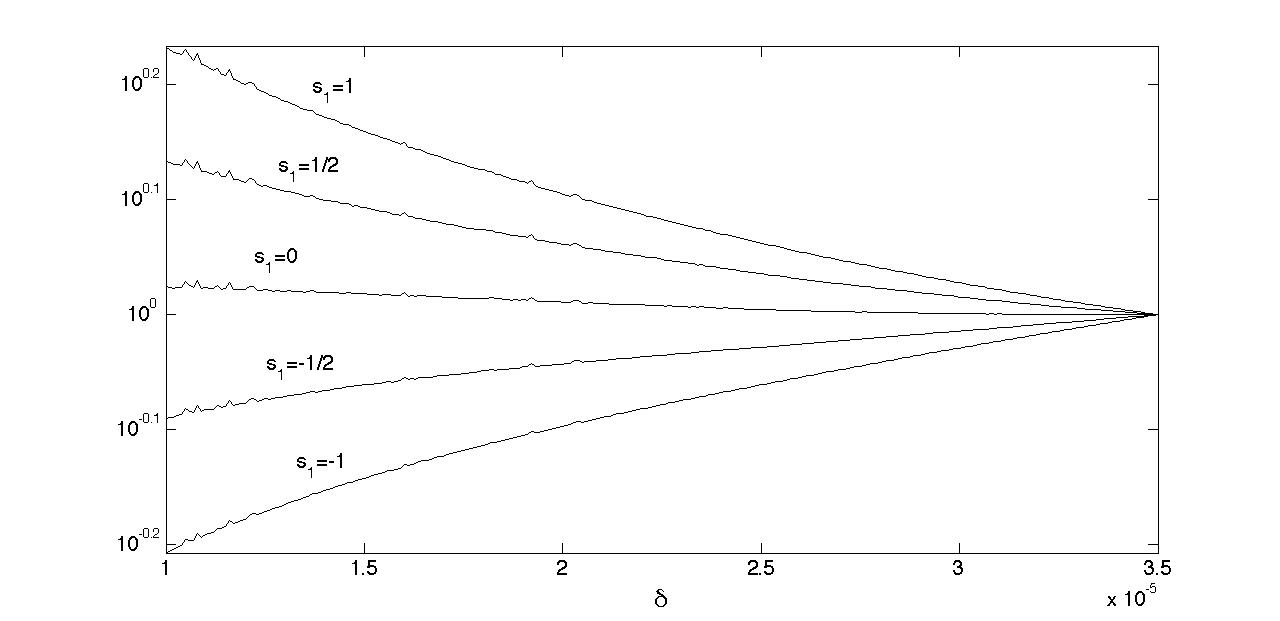}
\caption{\label{fig: norms}Normalised errors $c(s_1)\|u-u_\delta\|_{H^{s_1}(\T^1)}$ in logarithmic scale with different values of $s_1$. We observe that $u_\delta$ converge to $u$ at least when $s_1\leq -1/2$.}
\end{center}
\end{figure}

We know that $u, u_\delta\in H^1$ for all $\delta>0$ and are interested to know what happens in $H^1$ when $\delta\to 0$.  From figure  \ref{fig: norms} we can see that $u_\delta$ converges to $u$ in $H^{s_1}$ at least when $s_1<-1/2$. On the other hand even though both functions belong to  space $H^1$ we do not have convergence there that is $u_\delta\not\to u$ in $H^1$.

\section{Convergence of the finite-dimensional minimization problems}\label{sec:Gamma convergence}  

 {\newtext 
 In this section we consider the convergence of the finite-dimensional minimization problems (\ref{discrTikh1})
 to the ideal, infinite-dimensional model (\ref{contTikh0}).

Below, we consider
the case when $P_k$ is given by formula (\ref{Pk definition}) where
$\phi_j$ are the eigenfunctions of the Laplace operator of $N$ such that 
$(\phi_j)_{j=1}^\infty$ is an orthogonal basis  in  $L^2(N)$. Then $P_k:H^s(N)\to \R^k$ is a bounded linear map for all $s\in \R$. 

We will first consider the  relation of the noise models 
 (\ref{gauss}) and  (\ref{gauss infinite}).
Let $(\Omega,\Sigma,\mathbb P)$ be a probability space.
Assume that ${\mathcal E}:\Omega\to \mathcal D^\prime(Y)$ is a Gaussian random generalized
function (see \cite{B}) with a covariance operator $C_{\mathcal E}$ that 
can be extended to a bounded map $C_{\mathcal E}:H^{\widetilde{s}}(N)\to  H^{-\widetilde{s}}(N)$ with
some $\widetilde{s}\in \R$. Let us define $\widetilde{{\mathcal E}}= (I-\Delta)^{p/2+\widetilde{s}/2}{\mathcal E}$. Then the covariance operator of $\widetilde{{\mathcal E}}$ is $C_{\widetilde{{\mathcal E}}}=(I-\Delta)^{-p/2-\widetilde{s}/2}C_{\mathcal E}(I-\Delta)^{-p/2-\widetilde{s}/2}$. Weyl's theorem implies that  the eigenvalues of $-\Delta$ 
have the asymptotics $\lambda_j(-\Delta)= cj^{2/d}(1+O(j^{-1}))$ and thus $(I-\Delta)^{-p/2}$ is a Hilbert-Schmidt operator, that is $(I-\Delta)^{-p/2}\in \mathfrak{S}_2$, when $p>d/2$. On the other hand $(I-\Delta)^{-\widetilde{s}/2}C_{\mathcal E}(I-\Delta)^{-\widetilde{s}/2}: L^2 \to L^2$ since 
$C_{\mathcal E}: H^{\widetilde{s}}\to H^{-\widetilde{s}}$ is bounded and hence we can conclude that $C_{\widetilde{{\mathcal E}}}$ is a trace class operator which implies that ${\mathcal E}$ is almost surely a $H^{-p-\widetilde{s}}(N)$-valued random function.


Under the above assumptions, consider the $\R^k$-valued random variable
 $\tilde {\bf E}^{(k)}=P_k{\mathcal E}$. It has the  covariance operator $P_kC_{\mathcal E}P_k^*$. Thus 
the random variable
  $\tilde {\bf E}^{(k)}$  has the  distribution $N(0,\delta^2I)$ for all $k\in \Z_+$  if and only
  $C_W=\delta^2I$, or equivalently, ${\mathcal E}$ has the same distribution as $W\delta$, where $W:\Omega\to \mathcal D^\prime(Y)$ 
is the normalized Gaussian white noise given by (\ref{gauss infinite}).
This is the reason behind our assumption that in the ideal, infinite-dimensional 
model we have ${\rm noise}=\eps(x)\delta$ where $\eps$ is a realization of the normalized Gaussian white noise $W$.

To study the  convergence of the finite-dimensional minimization problems  (\ref{discrTikh1})
 to the infinite-dimensional problem (\ref{contTikh0}) we use the $\Gamma$-conver\-gence, see \cite{Attouch,Aubert,dalmaso}. Let $(Y,d_Y)$ be metric space and $\tau$ be the topology of $Y$ induced by the metric $d_Y$.
 Below we will in particular consider the case when $Y\subset X$
 is a closed bounded subset of  Banach space $X$ for which the dual space $X^\prime$ is separable.
 Then the weak topology of $X$ induces a topology $\tau$  for the subset $Y$
 that is metrizable, that is, induced by some metric $d_1:Y\times Y\to \R$, see
 see e.g. \cite[Prop. 8.7]{dalmaso}. 
Note that the metric $d_1$ is not necessarily the metric induced by the norm of $X$.
On following definition, see e.g. \cite[Def.\ 2.1.7]{Aubert} or \cite[Prop. 8.1]{dalmaso}.

\begin{definition}
  \label{def:gconv}
  Let $Y\subset X$ be a closed bounded subset of Banach space $X$ for which the dual space $X^\prime$ is separable.
  We say that $F_j:Y\to \R\cup\{\infty\}$ $\Gamma$-converges to $F:Y\to \R\cup\{\infty\}$ with the topology $\tau$
  and denote $F=\Glim_{j\to\infty} F_j$ if
  \begin{itemize}
  \item[(i)] For every $u\in Y$ and for every sequence $u_j$ $\tau$-converging to $u$ in $Y$ we have
      $F(u) \leq \liminf_{j\to\infty} F_j(u_j)$.
  \item[(ii)] For every $u\in Y$ there exists a sequence $u_j$ $\tau$-converging to $u$ in $Y$ such that
      $F(u) \geq \limsup_{j\to\infty} F_j(u_j)$.
  \end{itemize}
\end{definition}

We need also the concept of equicoercivity, see  \cite[Def.\ 2.1.8]{Aubert}.

\begin{definition}  
 Let $Y\subset X$ be a closed bounded subset of Banach space $X$ for which the dual space $X^\prime$ is separable.
  We call a sequence of functionals $F_j: Y\to \R\cup\{\infty\}$, $j\in\Z_+$, 
  {\em equicoercive} in topology $\tau$
  if for every $t\geq 0$ there exists
  a compact set $K_t\subset Y$ such that $\{u\in Y \; | \; F_j(u)\leq t\} \subset K_t$ for all $j\in\Z_+$.
\end{definition}

Using these definitions, we return to the  setting of the problem given in Section \ref{intro}
where
 $N$ is a $d$-dimensional compact closed manifold, $A\in\Psi^{-t}$ is a pseudodifferential operator, and
$m=m_\delta$ is the measurement (\ref{eq: measurement}).

Let us consider the
finite-dimensional minimization problems analogous to (\ref{discrTikh1}), that are given by 
 \begin{equation}\label{discrTikh1 modified}
  T_{\alpha;n,k}(m):=\argmin_{u\in X_n} F_{n,k}(u),
 \end{equation} 
  where $n,k\in \Z_+$ and $ F_{n,k}:H^r(N)\to \R\cup\{\infty\}$,
 \begin{equation}\label{discr funct. defined}  
    F_{n,k}(u)=\| P_kA u - P_km\|_{2}^2 +\alpha\|u\|_{H^r(N)}^2,\quad\hbox{for }u\in X_n, 
 \end{equation} 
 and  $F_{n,k}(u)=\infty $ for $u\not \in X_n$,
where $X_n\subset H^r(N)$ is a $n$-dimensional subspace. 
Also, let 
 $G:H^r(N)\to \R$ be   
 $$G(u)=\|A u\|_{L^2(N)}^2 - 2\langle m,A u \rangle +
 \alpha\|u\|_{H^r(N)}^2 ,\quad\hbox{for }u\in H^r(N).
 $$
 Let $Y=\{u\in H^r(N);\
\|u\|_{H^r(N)}\leq C_0\}$, where $$C_0>2\alpha^{-1} \max\bigg(\|A^*m\|_{H^{-r}(N)},\sup_{k\in \Z_+}\|(P_kA)^*m\|_{H^{-r}(N)}\bigg)$$ 
so that
$G(u)>G(0)=0$ and $F_{n,k}(u)>F_{n,k}(0)=\|P_k\varepsilon\|_{L^2}^2$  for all $\|u\|_{H^r(N)}> C_0$. 
Thus  functions $F_{n,k}:H^r(N)\to \R\cup\{\infty\}$ and $G:H^r(N)\to \R\cup\{\infty\}$
obtain their minimal values in $Y$.
We endow $Y$ with the relative topology determined by 
the weak topology of $H^{r}(N)$.

\begin{proposition}\label{prop 1}
Let the assumptions of Theorem \ref{speedelliptic} hold, in particular,
let $\varepsilon\in H^s(N)$ with $s<-d/2$ and $m=m_\delta$ be the measurement given by
 (\ref{eq: measurement}).
Moreover, 
assume 
that  $X_n\subset X_{n+1}$ and  $\cup_{n=1}^\infty X_n$ is a dense subset of $H^r(N)$  and let $c_k=\|P_km\|_{2}^2$.
  Then the functions $G_{n,k}:Y\to \R\cup\{\infty\}$, $$G_{n,k}(u)=F_{n,k}(u)-c_k,$$ converge to $G:Y\to \R$ as $n,k\to \infty$ in sense of  the $\Gamma$-convergence with respect to the topology of $Y$. Moreover,
  the minimizers  $   T_{\alpha;n,k}(m)$ of $F_{n,k}$ converge to the unique minimizer $T_{\alpha}(m)$
   of $G:H^r(N)\to \R$
 in the weak topology of $H^r(N)$ as $n,k\to \infty$.

\end{proposition}

\medskip

\noindent\textbf{Proof.} 
Let $u\in Y$ and let
 $u_{n,k}\in Y$ be a sequence that converge to $u$ weakly in $H^r(N)$ as $n,k\to \infty$.
As the linear operator $A:H^r(N)\to L^2(N)$ is a compact operator, $P_kAu_{n,k}$ converge 
to $Au$ in the strong topology of $L^2(N)$ as $n,k\to \infty$. Moreover,
 the map $u\mapsto \|u\|_{H^r(N)}$ is
 lower a semi-continuous function in $Y$. These facts imply that  
 the property (i) in Def.\ \ref{def:gconv} holds. 
 
Let  $Q_n$ be orthogonal projectors in  $H^r(N)$ onto the subspace $X_n$.  Let
 $u\in Y$, and define for $n,k\in \Z_+$
   $u_{n,k}=Q_nu$. Then 
   $$G_{n,k}(u_{n,k})=\| P_kA Q_nu - P_km\|_{2}^2 +\alpha\|Q_nu\|_{H^r(N)}^2-c_k$$ 
   converge to $G(u)$ as $n,k\to \infty$, and we see that
 the property (ii) in Def.\ \ref{def:gconv} is valid. Thus $G_{m,k}$ $\Gamma$-converge to $G$ as  
 $n,k\to \infty$.
 
%

Since all closed subsets of $Y$ are compact, we  see that $\{G_{n,k}:Y\to \R\cup\{\infty\};\ n,k\in \Z_+\}$  is an  equicoercive family of functions. Moreover,
the functions $G_{n,k}:Y\to \R\cup\{\infty\}$ and $F:Y\to \R$ have unique minimizers
and  the minimizer of $F_{n,k}:Y\to \R\cup\{\infty\}$ is equal to 
 the minimizers  $   T_{\alpha;n,k}(m)$ and finally, the minimizer of $G:Y\to \R$ is equal to $T_{\alpha}(m)$.
Thus by  \cite[Cor. 7.24] {dalmaso}, the minimizers  $   T_{\alpha;n,k}(m)$ of the functions $G_{n,k}$ converge weakly in $H^r(N)$ 
 to the minimizer of $G$ as $n,k\to\infty$.\hfill  $\square$
\medskip

\section{Conclusion}

\noindent We discuss above finite-dimensional linear models of indirect measurement corrupted by white Gaussian noise. Such models are used in countless practical inverse problems. It is desirable to connect these discrete models to an infinite-dimensional limit model. Such a connection can provide, for instance,  error analysis for numerical inversion and computational speed-ups based on robust switching between different discretizations related to multigrid methods.

The focus of our analysis is the apparent paradox arising from the (almost surely) infinite $L^2$-norm of the natural limit of white Gaussian noise in $\R^n$ as $n\rightarrow \infty$. We show how to build a rigorous theory removing this paradox, and we explain how to take this into account in discrete inverse problems using appropriate Sobolev space norms.

Proposition \ref{prop 1} shows that the infinite-dimensional minimization problem (\ref{contTikh1}) is the natural limit of the finite-dimensional minimization problems (\ref{discrTikh1 modified}). Therefore,
when the measured data is corrupted by white Gaussian noise,
despite the fact that the realizations of the white noise are almost surely not $L^2$ functions,
 the  inner product associated to the $L^2$-norm is appropriate for data fidelity terms when the
 inverse problems are solved
using Tikhonov regularization. Moreover, our results show how the
regularization parameters can be chosen to obtain converging results
when the noise amplitude goes to zero.

Our results pave the way to numerical analysis of Tikhonov regularization based on fruitful interplay between discrete and continuous models.
 
 }

\medskip
\textbf{Acknowledgements.}
This work was supported by the Finnish Centre of Excellence in Inverse Problems Research 2012-2017 (Academy of Finland CoE-project
250215). In addition, H.K. and M.L. were supported by Academy of Finland, project 
141104, and S.S. was supported by Academy of Finland, project 141094.

\bibliographystyle{amsalpha}

\vspace{1cm}

\end{document}